\documentclass[12pt]{article}
\usepackage{amsmath,amssymb,amsfonts,euscript,mathtools}
\usepackage{amsthm}
\usepackage{enumerate}
\usepackage{cite}
\usepackage[bookmarks]{hyperref}
\usepackage{cleveref}
\usepackage{libertine}
\usepackage[left=1.08in, right=1.08in]{geometry}
\usepackage{ifpdf}
\ifpdf    % we are running pdflatex
\usepackage{hyperref}
\else    % we are running latex
\usepackage[hypertex]{hyperref}
\fi

\newtheorem{prethm}{{\bf Theorem}}

\newtheorem{prelemma}[prethm]{Lemma}

\newcommand{\Gb}{{\overline G}}

\title{\large \bf Eccentricity and Algebraic Connectivity of Graphs}  

\author{{\normalsize
		{\sc B. Afshari},\,
		{\sc M. Afshari${}^{\mathsf{*}}$}
	}
	\vspace{3mm}
	\\{\footnotesize{${}^{\mathsf{*}}$\it Department of Computer Science and Engineering, Michigan State University}} {\footnotesize{}}
}
\date{}  
\begin{document}

\maketitle

\begin{abstract}
Let $G$ be a graph on $n$ nodes with algebraic connectivity $\lambda_{2}$. 
The eccentricity of a node is defined as the length of a longest shortest path starting at that node. 
If $s_\ell$ denotes the number of nodes of eccentricity at most $\ell$, then for $\ell \ge 2$,
	$$\lambda_{2} \ge \frac{ 4 \, s_\ell }{ (\ell-2+\frac{4}{n}) \, n^2 }.$$ 
As a corollary, if $d$ denotes the diameter of $G$, then
	$$\lambda_{2} \ge \frac{ 4 }{ (d-2+\frac{4}{n}) \, n }.$$ 

It is also shown that 
$$\lambda_{2} \ge \frac{ s_\ell }{ 1+ \ell \left(e(G^{\ell})-m\right) },$$
where $m$ and $e(G^\ell)$ denote the number of edges in $G$ and in the $\ell$-th power of $ G $, respectively. 
\end{abstract}

%\vspace{3mm}
\noindent {\em 2010 {\rm AMS} Classification}: 05C50, 15A18

%15A42
\bigskip
\bigskip
%\hspace{5mm}

Let $ G $ be a simple graph on $ V = \{v_1, v_2, \ldots , v_n\} $ with $m$ edges $E$. 
%We denote the \textit{degree} of $v_i$ by $d_i$. 
The \textit{distance} between nodes $v_i$ and $ v_j $, denoted by $ \mathrm{dist}(v_i,v_j) $, is the number of edges in a shortest path joining them. 
If there is no such path, then we define this value to be $\infty$. 
The \textit{diameter} of $G$, denoted by $d$ is the maximum distance between any pair of nodes of $G$. 
The $\ell$-th \textit{power} of $ G $, denoted by $ G^\ell $, is the graph with the same node set as $ G $ such that two nodes are adjacent in $ G^\ell $ if and only if their distance is at most $ \ell $ in $ G $. 
The \textit{eccentricity} of a node in $G$ is defined as the length of a longest shortest path starting at that node. 
We denote the number of nodes of eccentricity at most $\ell$ by $s_\ell$. 

The Laplacian matrix of $ G $ is defined as $L(G) = D-A$, where $ A $ is the adjacency matrix of $ G $ and $D$ is a diagonal degree matrix with $D(i,i) = \sum_{j} A(i,j)$. 
The smallest eigenvalue of $L(G)$ is 0, where the corresponding eigenvector is the all one vector. 
The second smallest eigenvalue of $L(G)$, denoted by $\lambda_{2}$, is called the \textit{algebraic connectivity} of $G$. 
	It is easy to see that the number of nodes with degree $n-1$ is a lower bound for this parameter, that is, 
		$$\lambda_{2} \ge s_1.$$
	
	Also it is known that \cite{AlgConSecondPow} 
	\begin{equation*}
		\lambda_{2} \ge \frac{ s_2 }{ n }.
	\end{equation*}

	Here it is shown that for $\ell \ge 3$, 
\begin{equation}\label{g1}
	\lambda_{2}
	\ge \frac{ 4 \, s_\ell }{ (\ell-2+\frac{4}{n}) n^2 }. 
\end{equation} 
This result can be viewed as a generalization, with a slight improvement, of the well-known bound~\cite{Mohar}: $$\lambda_2 \ge \frac{4}{d n},$$
as setting $\ell=d$ in (\ref{g1}) leads to $s_d = n$, and yields 
\begin{equation*}
	\lambda_2 \ge \frac{4}{(d-2+\frac{4}{n}) n}.
\end{equation*}

	It is also shown that 
	\begin{equation}\label{g2}
		\lambda_{2} 
		\ge \frac{ s_\ell }{ 1+ \ell \left(e(G^{\ell})-m\right) },
	\end{equation} 
	where $e(G^\ell)$ denotes the number of edges in $ G^\ell $. 
	This also generalizes the known result~\cite{Lu}:
	\begin{equation*}
		\lambda_2 \ge \frac{n}{1+ d \, e(\Gb)}.
	\end{equation*}

\section*{Proofs}
	\smallskip
	\begin{proof}[Proof of (\ref{g1})] 
If $G$ is not connected, then $s_\ell=0$, and the statement holds trivially. 
Therefore, we assume that $G$ in connected.  
		Let $\gamma = (\ell-2+\frac{4}{n})\frac{n^2}{4}$. 
		It is sufficient to show that $M= \gamma L(G)-L(G^\ell)$ is positive semi-definite. 
		Because if so, assuming that $z$ is the unit eigenvector corresponding to $\lambda_2$, we can deduce (\ref{g1}) as follows: 
		\begin{align*}
			\lambda_{2} 
			=z^T L(G) z 
			\ge \frac{1}{\gamma} z^T L(G^\ell) z
			\ge \frac{1}{\gamma} \lambda_{2}(G^\ell)
			\ge \frac{1}{\gamma} s_1(G^\ell) 
			= \frac{1}{\gamma} s_\ell.
		\end{align*}	
	
	Assume, to the contrary, that the minimum eigenvalue of $M$, denoted by $\mu$, is negative, and let $x = (x_1, x_2, \ldots , x_n)^T $ be an eigenvector corresponding to $\mu$. 
	Since the all one vector is an eigenvector corresponding to eigenvalue 0 of $M$, we have $\sum_{i=1}^n x_i = 0$. 
	We may assume that 
	$$ x_n \le x_{n-1} \le \ldots \le x_1. $$ 
	
		For $k\in \{1,2,\ldots,n-1\}$, we define $E_k$ as the set of edges $v_i v_j$ in $G$ where $i\le k<j$. 
		Let $k$ be chosen such that $$\sum_{v_i v_j\in E_k} |x_i-x_j|$$ is maximized, with the assumption that this sum equals 1 without loss of generality. 
		So, in particular, every edge $v_i v_j$ in $G$ satisfies $|x_i- x_j| \le 1$. %, as it belongs to some $E_q$.
		We proceed to show that
		\begin{equation}\label{e1}
			\sum_{ v_i v_j\in E(G^\ell):\, i\le k<j }  |x_i-x_j|  \le \gamma.
		\end{equation}
		
		Let $S= \{v_1,v_2,\ldots,v_k\}$. 
		For $r\ge 1$, we define $A_r$ as the set of nodes in $S$ whose shortest paths to $V\setminus S$ have length $r$. 
		Similarly, $B_r$ is the set of nodes in $V\setminus S$ whose shortest paths to $S$ have length $r$. 
		If $\mathrm{dist}(v_i,v_j)=q$ for some $i\le k<j$, then $v_i\in A_r$ and $v_j\in B_{r'}$ for some $r,r'\ge 1$ where $r+r'=q+1$. 
		Therefore, 
		\begin{equation*}
%			\sum_{i\le k<j:\, \mathrm{dist}(v_i,v_j)\le \ell} |x_i-x_j| 
			\sum_{ v_i v_j\in E(G^\ell):\, i\le k<j } |x_i-x_j| 
			\le \sum_{r,r' \ge 1: \, r+r'\le \ell+1} \sum_{v_i\in A_r, v_j\in B_{r'}} (x_i-x_j) .
		\end{equation*}

			\medskip
		Let $r,r'\ge 2$. 
		Let $s$ be the smallest number such that $v_{s}\in A_{r-1}$ and $t$ be the largest number such that $v_{t}\in B_{r'-1}$. 
We have
\smallskip
\begin{equation*}
	\sum_{ v_i\in A_r , v_j\in B_{r'} } ( x_i-x_j ) 
%		\sum_{ \substack{ v_i\in A_r \\ v_j\in B_{r'} } } ( x_i-x_j ) 
	\le \sum_{ v_i\in A_r , v_j\in B_{r'} } (x_{s}-x_{t}) + \sum_{v_i\in A_r}\sum_{ v_j\in B_{r'}:\, j> t } (x_{t}-x_j) + \sum_{v_j\in B_{r'}} \sum_{ v_i\in A_r:\, i< s } (x_i- x_{s}).
\end{equation*}
	
	\medskip
	There is a path $P$ of length $r-1$ from $v_{s}$ to a node in $V\setminus S$, and a path $P'$ of length $r'-1$ from $v_{t}$ to a node in $S$. 
	Both paths $P$ and $P'$ have an edge in $E_k$. 
	Therefore, we conclude:  
	\begin{align*}
		x_{s}-x_{t} 
		\le \sum_{v_i v_j \in E(P)\cup E(P')} |x_i-x_j| 
		\le r+r'-3 .
	\end{align*}

	Moreover, for each $v_i\in A_r$ with $i< s$, there is some edge $v_i v_j$ with $j\ge s$ in $G$. 
	Consequently, 
	\begin{align*}
		\sum_{v_i\in A_r: \, i < s} (x_i- x_{s}) 
		\le \sum_{v_i v_j\in E_{s-1}} |x_i-x_j| 
		\le 1.
	\end{align*}
Similarly, 
	\begin{align*}
		\sum_{v_j\in B_{r'}: \, j > t} (x_{t}-x_j) 
		\le \sum_{v_i v_j\in E_{t}} |x_i-x_j| 
		\le 1. 
	\end{align*} 
In consequence, 
\begin{align*}
	\sum_{ v_i\in A_r, v_j\in B_{r'} } (x_i-x_j) 
	&\le |A_r| |B_{r'}| (r+r'-3) + |A_r| +|B_{r'}| .
\end{align*}

\medskip
For each $v_j\in B_{1}$, let $v_{j'}$ be the node adjacent to $v_j$ in a shortest path from $v_{s}$ to $v_j$. 
Also, denote by $P_{sj'}$ a shortest path from $v_{s}$ to $v_{j'}$. 
We then have: 
\begin{align*}
	\sum_{ v_i\in A_r , v_j\in B_{1} } (x_i-x_j) 
	&\le \sum_{v_j\in B_{1}} \big(\sum_{ v_i\in A_r } (x_{s}-x_j) + \sum_{ v_i\in A_r:\, i< s } (x_i- x_{s})\big) \\
	&\le |A_r| \sum_{v_j\in B_{1}} \big((x_{s}-x_{j'})+(x_{j'}-x_j)\big) +|B_{1}|  \\
&\le |A_r| \sum_{v_j\in B_{1}} \sum_{v_a v_b\in P_{sj'}}|x_a-x_b| + |A_r| \sum_{v_i v_j\in E_k} |x_i-x_j| + |B_{1}|  \\
	&\le |A_r| |B_{1}| (r-2) + |A_r| +|B_{1}| .
\end{align*}
Similarly, 
\begin{align*}
	\sum_{ v_i\in A_1 , v_j\in B_{r'} } (x_i-x_j) 
	&\le |A_1| |B_{r'}| (r'-2) + |A_1| +|B_{r'}| .
\end{align*}

\medskip
Let $A_1 =\{v_{t_1},v_{t_2},\ldots,v_{t_p}\}$ and $B_1 =\{v_{t_{p+1}},v_{t_{p+2}},\ldots,v_{t_q}\}$ where $t_1 < t_2 <\ldots<t_q$. 
Note that for $i=1,2,\ldots,q-1$, we have 
$$|E_k \cap E_{t_i}| \ge \min\{i,q-i\},$$ 
as there are edges in $E_k$ with end nodes $v_{t_1},v_{t_2},\ldots,v_{t_i}$ if $i\le p$, and $v_{t_{i+1}},v_{t_{i+2}},\ldots,v_{t_{q}}$ if $i\ge p+1$.
Suppose the vector $y$ is obtained from $x$ by translation such that $y_{t_{ \lceil q/2 \rceil }}$ becomes zero. 
This leads to, 
\begin{align*}
	\sum_{ v_i\in A_1 , v_j\in B_1 } (x_i-x_j) 
	=\sum_{ v_i\in A_1 , v_j\in B_1 } (y_i-y_j) 
	\le \sum_{i<j} |y_{t_i}-y_{t_j}| 
	&= \frac{1}{2} \sum_{i=1}^{q}\sum_{j=1}^{q} |y_{t_i}-y_{t_j}| \\
	&\le \frac{1}{2} \sum_{i=1}^{q}\sum_{j=1}^{q} \big(|y_{t_i}|+|y_{t_j}|\big) \\
	&= q \sum_{i=1}^{q} |y_{t_i}| \\
	&= q \big(\sum_{i\le q/2}y_{t_i} - \sum_{i> q/2}y_{t_i}\big) \\
	&= q \sum_{i=1}^{q-1} \min\{i,q-i\} (y_{t_i}-y_{t_{i+1}}) \\
	&\le q \sum_{i=1}^{q-1} |E_k \cap E_{t_i}| (y_{t_i}-y_{t_{i+1}}) \\
	&= q \sum_{v_i v_j\in E_k} |y_i-y_j| 
	\le q.
\end{align*}

\medskip

Now, we can deduce that, 
\begin{equation*}
	\sum_{ v_i v_j\in E(G^\ell):\, i\le k<j } |x_i-x_j| 
	\le |A_1|+|B_1| + \sum_{r,r' \ge 1: \, 3\le r+r'\le \ell+1}{ (r+r'-3)|A_{r}| |B_{r'}|+ |A_{r}| + |B_{r'}| } .
\end{equation*}
The last term is bounded by the maximum possible value of 
$$a_1+b_1 + \sum_{r,r' \ge 1: \, 3\le r+r'\le \ell+1}{ (r+r'-3) a_r b_{r'} + a_r + b_{r'}}, $$
where $a_{r},b_{r}\ge 0$ for each $r$, subject to the constraints $\sum_{r=1}^{\ell}a_r= k$ and $\sum_{r'=1}^{\ell}b_{r'}= n-k$. 
By fixing the $b_{r'}$ values, this term reaches its maximum when $a_t=k$ for some $t$ and $a_r= 0$ for $r\neq t$. 
Subsequently, after fixing the values for $a_{r}$'s, the term is maximized when $b_{\ell-t+1}=n-k$ and $b_{r'}= 0$ for $r'\neq \ell-t+1$. 
Therefore, the inequality~(\ref{e1}) is derived as: 
\begin{align*}
	\sum_{ v_i v_j\in E(G^\ell):\, i\le k<j } |x_i-x_j|
	\le (\ell-2) k(n-k) + n
	\le (\ell-2) \frac{n^2}{4} + n 
	= \gamma .
\end{align*}

				Now, if $e_i \in \mathbb{R}^n$ is the $i$-th standard basis vector, then 
				\begin{align*}
					0 
					&\le  \gamma \sum_{v_i v_j\in E_k} |x_i-x_j| - \sum_{ v_i v_j\in E(G^\ell):\, i\le k<j } |x_i-x_j| \\
%					&= \gamma \sum_{v_i v_j\in E_k} |x_i-x_j| - \sum_{ v_i v_j\in E(G^\ell):\, i\le k<j } |x_i-x_j| \\
										&= \sum_{i\le k} \big( \gamma \sum_{v_i v_j\in E:\, j>k}{ (x_i - x_j) } - \sum_{v_i v_j\in E(G^\ell):\, j>k}{ (x_i - x_j) } \big)\\ 
					&= \sum_{i\le k} \big( \gamma \sum_{v_i v_j\in E}{ (x_i - x_j) } - \sum_{v_i v_j\in E(G^\ell)}{ (x_i - x_j) } \big) \\ 
					& = (e_1+e_2+\ldots+e_k)^T M x \\
					& = \mu \sum_{i=1}^{k} x_i \\
					&<0,
				\end{align*}
which is a contradiction. 
Therefore, the proof is now complete. 
\end{proof}

	\bigskip
	
		\begin{proof}[Proof of (\ref{g2})] 
			We need to show that $\left( 1+ \ell  \left(e(G^{\ell})-m\right) \right) L(G)-L(G^\ell)$ is positive semi-definite. 
		Consider a vector $x\in \mathbb{R}^n$, and let $v_r v_s \in E(G^\ell)$ be such that $|x_r - x_s|$ is maximum. 
		Denote by $P$ a shortest path between $v_r$ and $v_s$. 
		We then have: 
		\begin{align*}
			x^T L(G^\ell) x
			&= \sum_{v_i v_j\in E(G^\ell)} (x_i - x_j)^2 \\
			&\le x^T L(G) x + \left(e(G^{\ell})-m\right) (x_r-x_s)^2 \\
			&\le x^T L(G) x + \ell \left(e(G^{\ell})-m\right) \sum_{v_i v_j\in P} (x_i-x_j)^2 \\
			&\le \left( 1+ \ell \left(e(G^{\ell})-m\right) \right) x^T L(G) x. 
		\end{align*}
		This completes the proof. 	
	\end{proof}

\bibliographystyle{amsplain}

\end{document}